\documentclass[11pt]{article}
\usepackage{graphicx}
\usepackage{amssymb,amsmath}
\usepackage{pdfsync}
\usepackage[latin1]{inputenc}
\newtheorem{teorema}{Theorem}[section]
\newtheorem{prop}[teorema]{Proposition}
\newtheorem{coro}[teorema]{Corollary}

\newtheorem{remark}[teorema]{Remark}

\usepackage{makeidx,multicol}

\makeindex 

\def\R{\mathbb R}

\title{A lower bound for   the isoperimetric deficit}
\author{Julià Cufí and Agustí Reventós$^{1}$}
\date{}
\begin{document} \maketitle 



\abstract{In this paper we provide a Bonnesen-style inequality which gives a lower bound for the isoperimetric deficit corresponding to a closed convex curve in terms of some  geometrical invariants of this curve. Moreover we give a geometrical interpretation for the case when equality holds.}
\footnotetext[1]{Work partially supported by grants MTM2012-36378 and MTM2012-34834 (MEC).}
\section{Introduction}
Let $K$ be a plane compact convex set of area $F$ with boundary a curve $C=\partial K$ of length $L$. As it is well known, the isoperimetric inequality states $$F\leq \frac{L^{2}}{4\pi},$$ with equality only for discs.

Introducing the quantity $\Delta=L^{2}-4\pi F$, called the {\em isoperimetric deficit}, the above inequality can be written as  
 $\Delta \geq 0$. In some sense $\Delta$ 
measures the extend to which the convex set is away from a disc. It is interesting to know 
upper and lower bounds for $\Delta$ in terms of quantities associated to $K$.

\medskip

Hurwitz, in his paper about the use of Fourier series in some geometrical problems \cite{Hurwitz},  
proves the following inequality, which is a sort of reverse isoperimetric inequality and provides an upper bound for $\Delta$,
\begin{eqnarray}\label{hurt}0\leq \Delta \leq \pi |F_{e}|,\end{eqnarray}
where $F_{e}$ is the algebraic area enclosed by the evolute of $C$. Equality holds when $C$ is a circle or parallel to an astroid. 

Recall that the {\em evolute} of a plane curve is the locus of its centers of curvature or, equivalently,  the {\em envelope} of all the normals to this curve (i.e.,  the tangents to the evolute are the normals to the curve).  

\medskip

As for lower bounds, along the 1920's Bonnesen provided some inequalities of the type
$\Delta\geq B$, 
where $B$ is a non-negative quantity associated to the convex set vanishing only for circles. Moreover these quantities have a relevant geometrical meaning (see \cite{osserman}).

\medskip

In this note we prove a Bonnesen-style inequality which gives a lower bound for the isoperimetric deficit in terms of the difference between the area enclosed by the pedal curve of $C$ with respect to the Steiner point of $K$, and the area enclosed by $C$.


The {\em pedal curve} of a plane curve $C$ with respect to a fixed point $O$ is the locus of points $X$ so that the line $OX$ is perpendicular to the tangent to $C$ passing through $X$. The {\em Steiner point} of a plane convex set $K$, or the curvature centroid of $K$,  is the center of mass of $\partial K$ with respect to the density function that assigns to each point of $\partial K$ its curvature.

\medskip
Let $A$ be  the area enclosed by the pedal curve of $C=\partial K$ with respect to the Steiner point of $K$. In Theorem \ref{teorema31} it is proved that 
\begin{eqnarray}\label{ineq}\Delta\geq 3\pi(A-F).\end{eqnarray}
So, the quantity $3\pi(A-F)$ is a lower estimate of the isoperimetric deficit.  Since $A\geq F$,  this inequality  implies the isoperimetric one. Moreover $A=F$ only for circles and so $\Delta=0$ implies $C$ is a circle. 

We  point out that inequality \eqref{bnv} shows that Theorem \ref{teorema31} improves Theorem 4.3.1 in \cite{groemer}. Moreover our lower bound has a very clear geometric significance.  

\medskip

For the especial case of convex sets of constant width we obtain the inequality
\begin{eqnarray}\label{22set}\Delta\geq\frac{32}{9}\pi(A-F),\end{eqnarray}
which in turn improves inequality  in page 144 of \cite{groemer}, as  inequality \eqref{e17} shows.

\medskip

 We also consider when equality holds in \eqref{ineq} and  \eqref{22set} . In Corollaries \ref{coro33} and \ref{coro44} it is shown that this is so for circles and curves which are parallel to an astroid or to an hypocycloid of three cusps, respectively.

\medskip
Finally in Propositions \ref{prop12} and \ref{teorema52} we prove  that for convex curves  $C$ parallel to an astroid or an hypocycloid of three cusps,  the evolute of $C$ is similar, with ratio $2$ or $3$ respectively, to  
 the corresponding astroid or hypocycloid.
 
\section{Preliminaries}  
\subsection*{Support function}
A straight line $G$ in the plane is determined by the angle $\phi$ that the direction perpendicular to $G$ makes with the positive $x$-axis and the distance $p=p(\phi)$ of $G$ from the origin. The equation of $G$  then takes the form 
 
\begin{eqnarray}\label{01}x \cos\phi + y \sin \phi-p=0.\end{eqnarray}

Equation \eqref{01}, when $p=p(\phi)$ varies with $\phi$, is the equation of a family of lines. If we assume that the $2\pi$-periodic function  $p(\phi)$ is differentiable, the envelope of the family is obtained from \eqref{01} and the derivative of its left-hand side, as follows:

\begin{eqnarray}\label{02}-x \sin\phi + y \cos \phi-p'=0,\quad p'=dp/d\phi.\end{eqnarray}
From \eqref{01} and \eqref{02} we arrive at a  parametric representation of the envelope of the lines \eqref{01}:

\begin{eqnarray*}\label{03}x=p \cos\phi -p' \sin\phi,\quad\quad y=p\sin\phi+p'\cos\phi.\end{eqnarray*}
If the envelope is the boundary $\partial K$ of a convex set $K$ and the origin is an interior point of $K$, then $p(\phi)$ is called the {\em support function} of $K$ (or the support function of the convex curve $\partial K$).

Since $dx=-(p+p'')\sin\phi\,d\phi $ and $dy=(p+p'')\cos\phi\,d\phi$ (we here assume that the function $p$ is of class $C^2$),  arclength measure on $\partial K$ is given by

\begin{eqnarray}\label{04}ds=\sqrt{dx^2+dy^2}=|p+p''|\,d\phi\end{eqnarray}
and the radius of curvature  $\rho$ by \begin{eqnarray*}\label{05}\rho=\frac{ds}{d\phi}=|p+p''|.\end{eqnarray*}
It is well known (see for instance \cite{Santa}, page 3)  that a necessary and sufficient condition for a  periodic function $p$ to be the support function of a convex set $K$ is that $p+p''>0.$  Finally, it follows from \eqref{04} that the length of a closed convex curve that has support function $p$ of class $C^2$ is given by
\begin{eqnarray}\label{06}L=\int_{0}^{2\pi}p\,d\phi. \end{eqnarray}
The area of the convex set $K$ is expressed in terms of the support function by 

\begin{eqnarray}\label{07}F=\frac{1}{2}\int_{\partial K}pds=\frac{1}{2}\int_{0}^{2\pi}p (p+p'')\,d\phi=\frac{1}{2}\int_{0}^{2\pi}p^{2}\,d\phi -\frac{1}{2}\int_{0}^{2\pi}p'^{2}\,d\phi.\end{eqnarray}

\medskip
For any curve $C$ given by $(x(\phi),y(\phi))$, convex or not,  we will say that $p(\phi)$ is the {\em generalized support function}  of $C$ when \begin{eqnarray*}
x(\phi)&=&p(\phi)\cos(\phi)-p'(\phi)\sin(\phi),\\
y(\phi)&=&p(\phi)\sin(\phi)+p'(\phi)\cos(\phi).
\end{eqnarray*}
Note that $p(\phi)$ is not necessarily a distance, as it happens when we define the support function of a convex set.  In fact, $|p(\phi)|$ is the distance from the origin to the tangent to $C$ at the point $(x(\phi),y(\phi))$.

\medskip

It is easy to see that the  generalized support function $p_{e}(\phi)$ of the evolute of $C=\partial K$ is $p_{e}(\phi)=-p'(\phi+\pi/2)$, where $p(\phi)$ is the support function of $C$ (see \cite{ERev}). Hence, assuming $p(\phi)$ is a $C^{3}$-function, the algebraic area $F_{e}$ enclosed by the evolute of $C$ is given by
$$F_{e}=\frac{1}{2}\int_{0}^{2\pi}p'(p'+p''')\,d\phi=\frac{1}{2}\int_{0}^{2\pi}p'^{2}\,d\phi-\frac{1}{2}\int_{0}^{2\pi}p''^{2}\,d\phi.$$

\subsection*{Steiner point}
The {\em Steiner point} $S(K)$ of a convex set $K$ of the Euclidean plane is defined by
$$S(K)=\frac{1}{\pi}\int_{0}^{2\pi}(\cos \phi,\sin\phi)p(\phi)d\phi,$$ where $p(\phi)$ is the support function of $\partial K$ (see \cite{groemer}).

Thus, if
\begin{eqnarray}\label{fourier}p(\phi)=a_{0}+\sum_{n\geq 1}a_{n}\cos n\phi+b_{n}\sin n\phi,\end{eqnarray}
is the Fourier series  of the $2\pi$-periodic function $p(\phi)$, 
the Steiner point is $$S(K)=(a_{1},b_{1}).$$

The Stiener point of $K$ is also known as the curvature centroid of $K$ because under appropriate smothness conditions it is the center of mass of $\partial K$ with respect to the density function that assigns to each point of $\partial K$ its curvature.

\medskip
The relation between the support function $p(\phi)$ of a convex set $K$ and the support function $q(\phi)$ of the same convex set  but with respect to a new reference with origin at the point $(a,b)$, and axes parallel to the previous 
$x$ and $y$-axes,  is given by 

$$q(\phi)=p(\phi)-a\cos\phi-b\sin\phi.$$ 

Hence, taking the Steiner point as a new origin, we have
$$q(\phi)=a_{0}+\sum_{n\geq 2}a_{n}\cos n\phi+b_{n}\sin n\phi.$$

We recall that the {\em Steiner ball} of $K$ is the ball whose center is the Steiner point and whose diameter is the mean width of $K$.

\subsection*{Pedal curve}

If the curve $C$ is given in cartesian coordinates as the envelope of the lines $x\cos\phi+y\sin\phi-p(\phi)=0$, then the pedal curve ${\cal P}={\cal P}(\phi)$ of $C$ with respect to the origin, is given by

 $${\cal P}(\phi)=(p(\phi)\cos\phi,p(\phi)\sin\phi),$$
or, in polar coordinates, by $r=p(\phi)$.

In particular, if $C$ is closed, the area enclosed by ${\cal P}$  
is \begin{eqnarray}\label{qf}A=\frac{1}{2}\int_{0}^{2\pi}p^{2}\,d\phi.\end{eqnarray}
If $F$ is the area enclosed by $C$, we obviously have  $A\geq F$  with equality if and only if $C$ is a circle.

\section{A lower bound for the  isoperimetric deficit}
We proceed now to provide a lower bound for the isoperimetric deficit.   
\begin{teorema}\label{teorema31}
Let $K$ be a compact convex set of area $F$ with boundary a curve $C=\partial K$ of class ${\cal C}^{2}$ and length $L$. Let $A$ be the area enclosed by the pedal curve of $C$ with respect to the Steiner point $S(K)$. Then
\begin{eqnarray}\label{poi}\Delta\geq 3\pi(A-F),\end{eqnarray} 
where $\Delta =L^2-4\pi F$ is the isoperimetric deficit.
\end{teorema}
{\em Proof}.  Let $p(\phi)$ be the support function of $C$,  with respect to an orthonormal reference with origin in the Steiner point,  and axes parallel to the $x$ and $y$-axes.

We know  that the Fourier series  of  $p(\phi)$, is 
$$p(\phi)=a_{0}+\sum_{n\geq 2}a_{n}\cos n\phi+b_{n}\sin n\phi.$$

By Parseval's identity we have
\begin{eqnarray}\label{parseval}\frac{1}{2\pi}\int_{0}^{2\pi}p^{2}\,d\phi=a_{0}^{2}+\frac{1}{2}\sum_{n\geq 2}(a_{n}^{2}+b_{n}^{2}),\end{eqnarray}
and similar expressions for $p'$ and $p''$. Concretely we have
\begin{eqnarray}\label{eee}\int_{0}^{2\pi}p'^{2}\,d\phi=\pi\sum_{n\geq 2}n^{2}(a_{n}^{2}+b_{n}^{2}), \qquad \int_{0}^{2\pi}p''^{2}\,d\phi= \pi\sum_{n\geq 2}n^{4}(a_{n}^{2}+b_{n}^{2}).\end{eqnarray}

Hence, the isoperimetric deficit $\Delta=L^{2}-4\pi F$, according to \eqref{06} and \eqref{07}, is given by

\begin{eqnarray*}\label{221}\Delta&=&\left(\int_{0}^{2\pi}p\,d\phi\right)^{2}-2\pi\int_{0}^{2\pi}p^{2}\,d\phi+2\pi\int_{0}^{2\pi}p'^{2}\,d\phi\\&=&2\pi^{2}\sum_{n\geq 2}(n^{2}-1)(a_{n}^{2}+b_{n}^{2})\geq \frac{3\pi^{2}}{2}\sum_{n\geq 2}n^{2}(a_{n}^{2}+b_{n}^{2})=\frac{3\pi}{2}\int_{0}^{2\pi}p'^{2}\,d\phi.\end{eqnarray*}

But it follows from \eqref{07} and \eqref{qf} that
$$\frac{1}{2}\int_{0}^{2\pi}p'^{2}\,d\phi=\frac{1}{2}\int_{0}^{2\pi}p^{2}\,d\phi- F=A-F,$$
and hence
\begin{eqnarray*}\Delta\geq 3\pi(A-F).\qquad \square\end{eqnarray*}

\medskip

The above proof shows that $\Delta=0$ if and only if $p(\phi)=a_{0}$, that is, when $C$ is a circle.

\bigskip

Although the constant $3\pi$ appearing in Theorem \ref{teorema31} cannot be improved for general convex sets, it is possible to obtain a stronger inequality
 for special type of convex sets.

 For instance,  for convex sets of constant width we have
 the following result.
 
 \begin{prop}\label{prop32}
 Let $K$ be a compact convex set of constant width in the hypothesis of Theorem \ref{teorema31}. Then
 $$\Delta\geq\frac{32}{9}\pi(A-F).$$
 \end{prop}
{\em Proof}.  Since constant width means 
$p(\phi)+p(\phi+\pi)$ is constant and 
$$p(\phi)+p(\phi+\pi)=2\sum_{0}^{\infty}(a_{2n}\cos 2n\phi+b_{2n}\sin 2n\phi)$$ 
it follows that $a_{n}=b_{n}=0$ for all even $n>0$.

Introducing this in the proof of Theorem \ref{teorema31} the result follows. $\square$

\subsection*{Relationship with the $L^{2}$ metric}

Consider now the quantity $\delta_{2}(K)$ equal to the distance in $L^{2}(S^{1})$ between the support function of $K$ and the support function of the Steiner ball of $K$.

It is known, see \cite{groemer} Theorem 4.3.1,
 that
\begin{eqnarray}\label{ghj}\Delta\geq 6\pi\delta_{2}(K)^{2}.\end{eqnarray}

We can state now the following inequality \begin{eqnarray}\label{bnv}3\pi(A-F)\geq 6\pi\delta_{2}(K)^{2}.\end{eqnarray}
 
To proof this we first observe that $$\delta_{2}(K)^{2}=\pi\sum_{n=2}^{\infty}(a_{n}^{2}+b_{n}^{2}),$$ where $a_{n},b_{n}$ are Fourier coefficients of the support function of $K$, $p(\phi)$, as in \eqref{fourier} (see page 142 of \cite{groemer}). 
 Moreover the proof of Theorem \ref{teorema31} shows that
 $$3\pi(A-F)=\frac{3\pi^{2}}{2}\sum_{n\geq 2}n^{2}(a_{n}^{2}+b_{n}^{2}),$$
 hence inequality \eqref{bnv}
follows.  
 
 So we have
 $$\Delta\geq3\pi(A-F)\geq 6\pi\delta_{2}(K)^{2}$$ 
 which improves the inequality  \eqref{ghj}.

 \medskip 
  For compact convex sets of constant width it is known that
  
  \begin{eqnarray}\label{ghjk}\Delta\geq 16\pi\delta_{2}(K)^{2}\end{eqnarray}
 (see page 144 of \cite{groemer}).
 
 We can state now the following inequality 
 
\begin{eqnarray}\label{e17}\frac{32}{9}\pi(A-F)\geq 16\pi\delta_{2}(K)^{2}.\end{eqnarray}
 The proof 
 is the same as for \eqref{bnv} taking into account that now the even coefficients vanish.

  So we have
 $$\Delta\geq \frac{32}{9}\pi(A-F)\geq 16\pi\delta_{2}(K)^{2}$$ 
 which improves the inequality  \eqref{ghjk}.

\section{Equality of the lower bound with the isoperimetric deficit} 

Now we study  the case of equality in Theorem \ref{teorema31}.
It is clear from the proof that $\Delta=3\pi(A-F)$ if and only if  $$p(\phi)=a_{0}+a_{2}\cos 2\phi+b_{2}\sin 2 \phi.$$
In order to characterize the curves with this type of support function we recall that  the parametric equations of an {\em astroid} (a $4$-cusped hypocycloid) are
\begin{eqnarray*}
x(\phi)&=&2a\sin^{3}(\phi),\\
y(\phi)&=&2a\cos^{3}(\phi),
\end{eqnarray*}for some constant $a\in\R$, $a\neq 0$, with $0\leq \phi\leq 2\pi$.
From this it is easy to see that the generalized support function $p(\phi)$ of the astroid  is $p(\phi)=a\sin(2\phi)$, 
where $\phi$ is the angle between the normal $(-y'(\phi),x'(\phi))$ and the positive $x$-axis.

This implies that the curves with generalized support function given by $$q(\phi)=b+p(\phi)=b+a\sin(2\phi),$$
where $b\in\R$, are parallel to an astroid. The distance between these curves and the astroid is $|b|$.

We have the following result.
\begin{prop}\label{prop1} Let $$p(\phi)=a_{0}+a_{2}\cos(2\phi)+b_{2}\sin(2\phi)$$ be the  support function of a closed convex curve $C$ of length $L$, with $a_{2}^2+b_{2}^2\neq 0$. Then
the interior parallel curve to $C$ at distance $L/2\pi$ is an astroid.
\end{prop}
{\em Proof.}
We make the change of variable $u=\phi-\phi_{0}+\frac{\pi}{4}$, where $$\tan 2 \phi_{0}=\frac{b_{2}}{a_{2}}.$$

Then

\begin{eqnarray*}\sin 2u&=&\cos 2(\phi-\phi_{0})=\cos 2\phi\cos 2\phi_{0}+\sin 2\phi\sin 2\phi_{0}\\&=&\cos 2\phi\frac{a_{2}}{\pm\sqrt{a_{2}^{2}+b_{2}^{2}}}+\sin 2\phi\frac{b_{2}}{\pm\sqrt{a_{2}^{2}+b_{2}^{2}}}\end{eqnarray*}

Hence

$$p(u)=a_{0}\pm a\sin 2u$$
where $a=\sqrt{a_{2}^{2}+b_{2}^{2}}$. This shows that the given curve is parallel to an astroid at distance $|a_{0}|$. By the condition of convexity, $p+p''=a_{0}\mp 3a\sin 2u >0$, and so $a_{0}$  is positive.  Since $L=\int_{0}^{2\pi}p(\phi)\,d\phi=2\pi a_{0}$,  the proposition is proved. $\square$


\begin{coro}\label{coro33}
Equality in Theorem \ref{teorema31} holds if and only if $C$ is a circle or a curve parallel to an astroid.
\end{coro}
{\em Proof}. As we have said, equality in \eqref{poi} holds when 
$$p(\phi)=a_{0}+a_{2}\cos 2\phi+b_{2}\sin 2 \phi.$$
If $a_{2}=b_{2}=0$, $p(\phi)=a_{0}$ is the support function of a circle. If $a_{2}^2+b_{2}^2\neq 0$, the result follows from Proposition \ref{prop1}. $\square$ 

\bigskip

Now we study the case of equality in Proposition \ref{prop32}.

It is clear from the proof of this Proposition  
 that equality holds if and only if $$p(\phi)=a_{0}+a_{3}\cos 3\phi+b_{3}\sin 3\phi.$$

 In order to characterize the curves with this type of support function we recall  that the parametric equations of an {\em hypocycloid of three cusps}, with respect to a  suitable  orthogonal system, are   
\begin{eqnarray*}
x(t)&=&-2a\cos t-a\cos 2t
\\
y(t)&=&-2a\sin t+a\sin 2t
\end{eqnarray*}
with $a\in \R$, $a\neq 0$, $t\in[0,2\pi]$.

%

The relationship between the parameter $t$ and the angle $\phi(t)$  between  the normal vector $(-y'(t),x'(t))$ and the positive $x$-axis
is
$$\phi(t)=\alpha(t)-\frac{\pi}{2},$$
where $\alpha(t)$ denotes the angle between the tangent vector $(x'(t),y't)$  and the positive $x$-axis.

Hence

$$\tan\phi(t)=-\cot\alpha(t)=\frac{\sin t+\sin 2t}{\cos t-\cos 2t}=\cot\frac{t}{2}$$
and so 

$$t=\pi-2\phi(t).$$


On the other hand, the generalized support function  $p(\phi)$ of the hypocycloid must verify

$$\left(\begin{array}{c}x(\phi) \\y(\phi)\end{array}\right)=\left(\begin{array}{cc}\cos\phi & -\sin\phi \\\sin\phi & \cos\phi\end{array}\right)\left(\begin{array}{c}p(\phi) \\p'(\phi)\end{array}\right),
$$
so

$$\left(\begin{array}{c}p(\phi) \\p'(\phi)\end{array}\right)=\left(\begin{array}{cc}\cos\phi & \sin\phi \\-\sin\phi & \cos\phi\end{array}\right)\left(\begin{array}{c}-2a\cos(\pi-2\phi)-a\cos2(\pi-2\phi) \\-2a\sin (\pi-2\phi)+a\sin2(\pi-2\phi)\end{array}\right).$$

Then, using standard addition trigonometric formulas, it follows

$$p(\phi)=a\cos(3\phi).$$

\begin{prop}\label{prop52}Let $$p(\phi)=a_{3}\cos3\phi+b_{3}\sin3\phi$$ be the generalized support function of a closed curve $C$, with $a_{3}^{2}+b_{3}^{2}\neq 0$.
Then $C$ is a hypocycloid  of three cusps. 
\end{prop}
{\em Proof}. We make the change of variable given by $u=\phi-\phi_{0}$,  where $$\tan 3\phi_{0} =\frac{b_{3}}{a_{3}}.$$

Then, an easy computation gives  $$p(u)=a\cos(3u),$$
where $a=\dfrac{a_{3}}{\cos\phi_{0}}$,  and the proposition follows. $\square$
\begin{coro}\label{coro44}
Equality in Proposition \ref{prop32} holds if and only if $C$ is a circle or a curve parallel 
to an hypocycloid of three cusps.\end{coro}
{\em Proof.} We have seen that equality in Proposition \ref{prop32} holds when 
$$p(\phi)=a_{0}+a_{3}\cos3\phi+b_{3}\sin3\phi.$$
If $a_{3}=b_{3}=0$, $p(\phi)=a_{0}$ is the support function of a circle.
If $a_{3}^{2}+b_{3}^{2}\neq 0$, the result follows from Proposition \ref{prop52}. $\square$

\begin{remark}\label{remarc34}\em 
As it is well known (see for instance page 8 of \cite{Santa})  the area $F_{r}$ enclosed by the interior parallel at distance $r$ to a closed curve is given by
\begin{eqnarray*}
F_{r}=F-Lr+\pi r^{2}
\end{eqnarray*}
 where $L$ and $F$
are respectively the length and the area corresponding to the given curve.

In particular, if $r=L/2\pi$, we get $$F_{L/2\pi}=F-\frac{L^{2}}{4\pi},$$ 
or, equivalently 
$$ \Delta=- 4\pi F_{L/2\pi},$$
which gives a geometrical interpretation of the isoperimetric deficit.

In particular   the isoperimetric inequality $\Delta\geq 0$ is equivalent to $F_{L/2\pi}\leq 0$, a fact that suggests  a more geometric proof of the isoperimetric inequality, by showing that in the process of collapsing,  the curve reverses orientation.  Moreover, $F_{L/2\pi}=0$ holds only for   circles. 
\end{remark}

\begin{remark}\em 
 Combining Theorem \ref{teorema31} with Hurwitz's inequality \eqref{hurt} we have
 the relation 
 \begin{eqnarray}\label{yyy}A-F\leq\frac{1}{3}|F_{e}|,\end{eqnarray}
with equality for circles or curves parallel to an astroid.
\end{remark}

%
%
%
%

\section{Parallel curves and evolutes}
We have seen the role played by the convex curves parallel to an astroid or to an hypocycloid of three cusps.  For such curves $C$ we show that there is a quite surprising relationship between the parallel curve at distance $L/2\pi$ to $C$ and the evolute of $C$.

\begin{prop}\label{prop12} Let $$p(\phi)=a_{0}+a_{2}\cos(2\phi)+b_{2}\sin(2\phi)$$ be the  support function of a closed convex curve $C$ of length $L$. Then the evolute  of $C$ and the interior  parallel curve to $C$ at distance $L/2\pi$,  are similar with ratio $2$.     
\end{prop}
{\em Proof.} We shall see that there is a similarity, composition of a  rotation with an homothecy, applying the parallel curve on the evolute.
We may assume, by the proof of Proposition \ref{prop1}, $p(\phi)=a_{0}+a\sin(2\phi)$.  
The generalized support function of the parallel curve to $C$ at distance $L/2\pi=a_{0}$ is $q(\phi)=a\sin(2\phi)$ and the corresponding one to the  evolute of $C$ is $$p_{e}(\phi)=-p'(\phi+\frac{\pi}{2})=2a\cos(2\phi).$$

The generalized support function of the rotated $3\pi/4$ parallel curve
is
$$\tilde{p}(\phi)=q(\phi-\frac{3\pi}{4})=a\cos(2\phi).$$

Hence this rotated curve is homothetic, with ratio $2$,  to the evolute. $\square$ 

\begin{remark}\em 
In particular, the area of the evolute of such a curve is four times the area of the parallel curve at distance $L/2\pi$. The reciprocal is also true.  In fact,  since  Hurwitz's inequality, by Remark \ref{remarc34},  is equivalent to 
$$4|F_{L/2\pi}|-|F_{e}|\leq 0,$$
the curves for which 
the area of the evolute is four times the area of the parallel curve at distance $L/2\pi$, are {\em exactly} circles or  curves parallel to an astroid.  
\end{remark}

\begin{prop}\label{teorema52} Let $$p(\phi)=a_{0}+a_{3}\cos(3\phi)+b_{3}\sin(3\phi)$$ be the  support function of a closed convex curve $C$ of length $L$. Then the evolute  of $C$ and the interior  parallel curve to $C$ at distance $L/2\pi$,  are similar with ratio $3$.     
\end{prop}
{\em Proof.} Analogous to that of Proposition \ref{prop12} with $q(\phi)=a\cos (3\phi)$, 
according to the proof of Proposition \ref{prop52}. $\square$

\medskip 
Next figures  show convex curves with support functions $p(\phi)=5+\sin(2\phi)$  and $p(\phi)=8+\sin(3\phi)$, their   parallel interior curves at  distance $L/2\pi=5$ and $L/2\pi=8$ respectively, and the corresponding evolutes.

\begin{center}
\includegraphics[width=.30\textwidth]{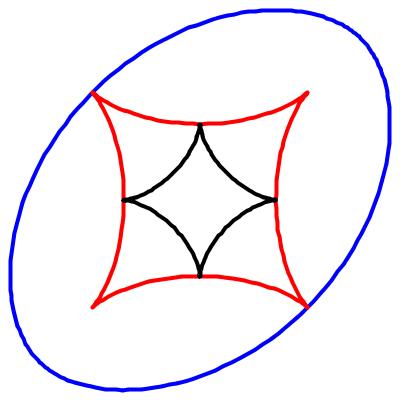}
\hspace{2cm}\includegraphics[width=.32\textwidth]{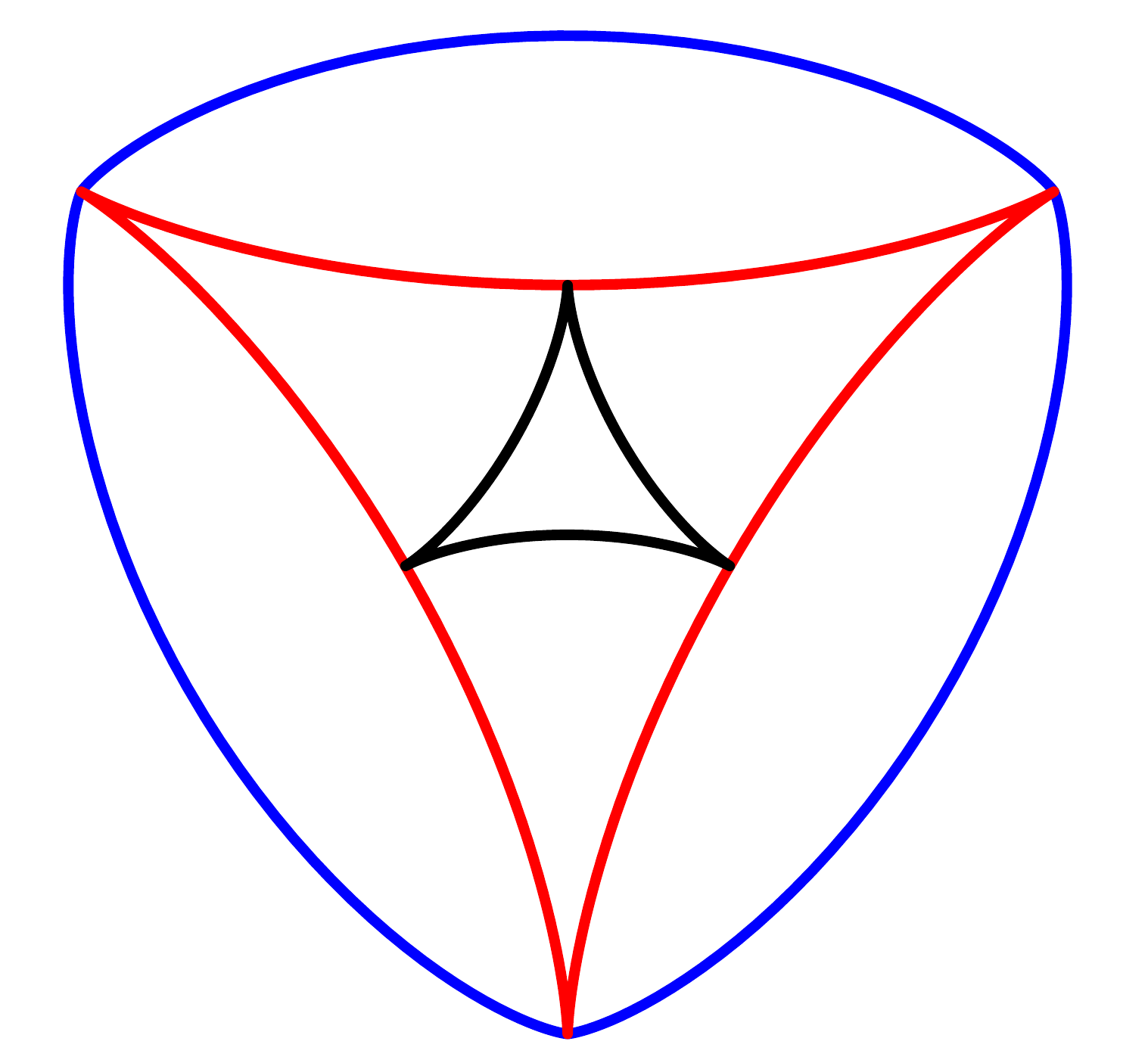}
\end{center}
\centerline{$p(\phi)=5+\sin(2\phi)$\hspace{4.2cm}$p(\phi)=8+\sin(3\phi)$}

\bibliographystyle{plain}
\bibliography{BibliotecaGNE}

\noindent {\em Departament de Matemàtiques \\Universitat Aut\`{o}noma de Barcelona\\ 08193 Bellaterra, Barcelona\\Catalonia\\

\noindent jcufi@mat.uab.cat, agusti@mat.uab.cat.

}

 \end{document}